\documentclass[reqno,11pt]{amsart}

\usepackage{epsf}
\usepackage{amssymb}
\usepackage{amsmath}

\date{}

\def\pic#1#2{ \def\epsfsize##1##2{1.0##1}\raisebox{#1}{$\vcenter{\hbox{\epsffile{#2.eps}}}$}}

\def\smallpic#1#2{ \def\epsfsize##1##2{0.9##1}\raisebox{#1}{$\vcenter{\hbox{\epsffile{#2.eps}}}$}}

\def\tablepic#1#2{ \def\epsfsize##1##2{0.2##1}\raisebox{#1}{$\vcenter{\hbox{\epsffile{#2.eps}}}$}}

\theoremstyle{plain}
\newtheorem{theorem}{Theorem}

\theoremstyle{definition}

\theoremstyle{remark}

\newtheorem*{remarks}{Remarks}

\def\C{{\mathbb C}}

\author{Sebastian Baader}
\title{A Note on Vassiliev invariants of quasipositive knots}

\begin{document}

\begin{abstract} It has been known that any Alexander polynomial of a knot can
be realized by a quasipositive knot. As a consequence, the Alexander
polynomial cannot detect quasipositivity. In this paper we prove a similar
result about Vassiliev invariants: for any oriented knot $K$ and any natural
number $n$ there exists a quasipositive knot $Q$ whose Vassiliev invariants of
order less than or equal to $n$ coincide with those of $K$.
\end{abstract}

\maketitle

A quasipositive braid is a product of conjugates of a positive standard
generator of the braid group. If a link can be realized as the closure of a
quasipositive braid then we call it quasipositive. When Lee Rudolph
introduced quasipositive links (in \cite{Ru2}), he showed that they could be
realized as transverse $\C$-links, i.e. as transverse intersections of complex
plane curves with the standard sphere $S^3 \subset \C^2$. Here a complex plane
curve is any set $f^{-1}(0) \subset \C^2$, where $f(z,w) \in \C[z,w]$ is a
non-constant polynomial. Conversely, every transverse $\C$-link is a
quasipositive link, as was recently proved by M.~Boileau and S.~Orevkov (in
\cite{BO}). We shall use yet another description of quasipositive knots which
is based upon Seifert diagrams in order to prove the following result.

\begin{theorem} For any oriented knot $K$ and any natural number $n$ there
exists a quasipositive knot $Q$ whose Vassiliev invariants of order
less than or equal to $n$ coincide with those of $K$.
\end{theorem}

This is related to a result of Lee Rudolph (\cite{Ru1}),
who showed that any Alexander polynomial can be realized by a quasipositive
knot. We also mention that theorem~1 was formulated as a question by
A.~Stoimenow in~\cite{St}.

The proof of theorem~1 is based upon a construction of Y.~Ohyama, who showed
that any finite number of Vassiliev invariants can be realized by an
unknotting number one knot (see \cite{O}). His construction involves certain
$C_n$-moves, which were defined by K.~Habiro in~\cite{Ha}, see
also~\cite{OTY}. A special $C_n$-move is defined diagrammatically in figure~1.
It takes place in a section with $2(n+1)$ endpoints or $(n+1)$ strands,
respectively. The strands are numbered from $1$ to $n+1$ and are all connected
outside the indicated section, since they belong to one knot $K$.
Going along $K$ according to its orientation, starting at the first strand, we
encounter the other strands in a certain order which depends on how the
strands are connected outside the indicated section. This order defines a
permutation, say $\sigma \in S_{n}$, of the numbers 2, 3, $\ldots$, n+1.

\begin{figure}[ht]
\pic{-0pt}{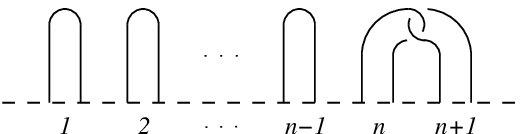} $\longleftrightarrow$ \pic{-0pt}{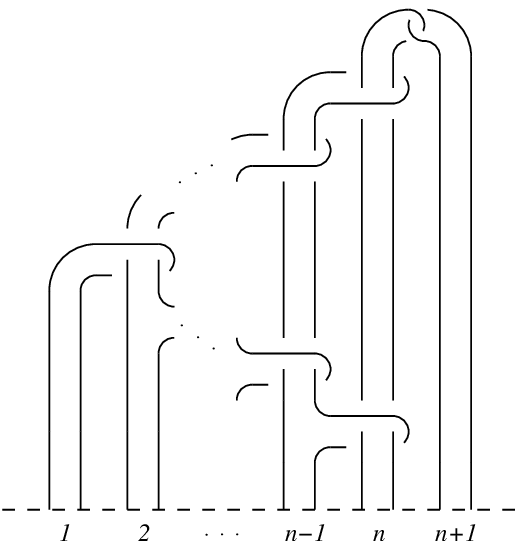}
\caption{}
\end{figure}

In~\cite{OT}, Y.~Ohyama and T.~Tsukamoto explain the effect of a $C_n$-move on
Vassiliev invariants of order $n$. Their result (\cite{OT}, theorem~1.2)
implies the following:

\begin{enumerate}
\item A $C_n$-move does not change the values of Vassiliev invariants of order
less than $n$.\\
\item Let $K$ and $\widetilde K$ be two knots which differ by one $C_n$-move,
and $v_n$ any Vassiliev invariant of order $n$. Then $|v_n(K)-v_n(\widetilde
K)|$ depends only on the permutation $\sigma \in S_{n}$ defined by the cyclic
order of the $(n+1)$ strands of the section where the $C_n$-move takes place.
\end{enumerate}

\begin{proof}[Proof of theorem~1]
Starting from the diagram of the positive twist knot $5_2$ shown in figure~2,
we construct a quasipositive knot $Q$ with the desired properties by applying
several $C_n$-moves, $2 \leqslant i \leqslant n$, step by step.

\begin{figure}[ht]
\pic{-0pt}{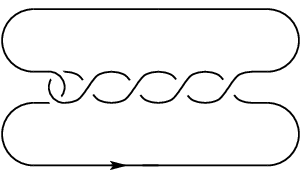} 
\caption{}
\end{figure}

In the first step, we construct a quasipositive knot $Q_1$ whose Vassiliev
invariants of order two (the Casson invariant) equals that of $K$. Choose
natural numbers $a$ and $b$, such that $v_2(K)=2+a-b$. Here $v_2(K)$ is the
Casson invariant of $K$. Using these two numbers, we define a knot $Q_1$
diagrammatically, as shown in figure~3. 

\begin{figure}[ht]
\pic{-0pt}{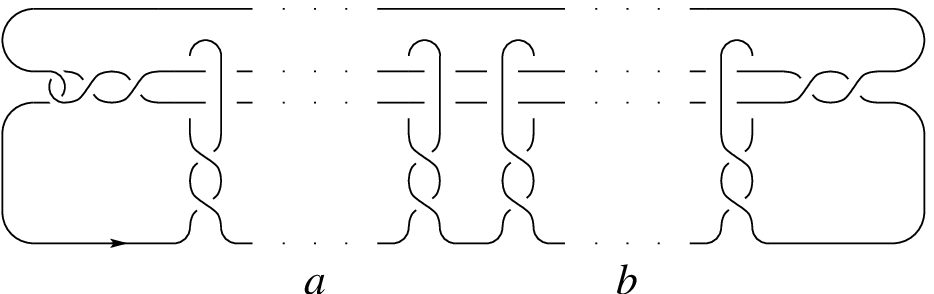}

\bigskip
\pic{-0pt}{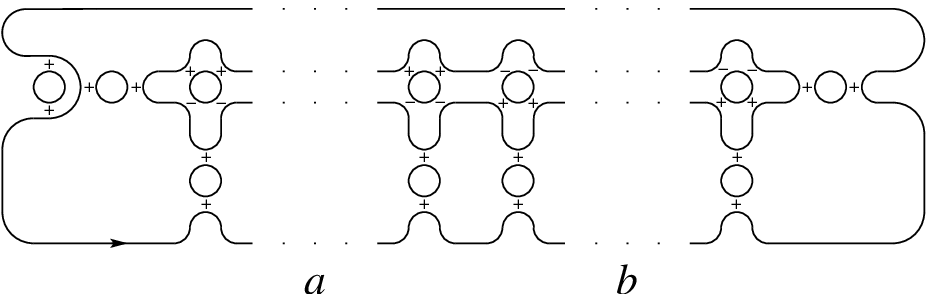}
\caption{}
\end{figure}

By construction, we have
$$v_2(Q_1)=2+a-b=v_2(K).$$
This follows easily by one application of the following relation for the
Casson invariant of knots:
$$v_2(\tablepic{-0pt}{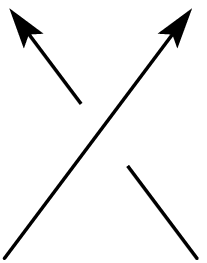})-v_2(\tablepic{-0pt}{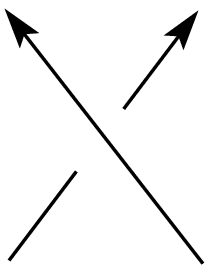})=lk(\tablepic{-0pt}{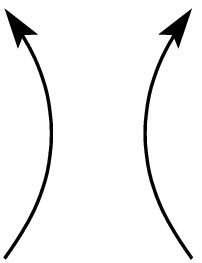}).$$

Indeed, a crossing change at the clasp on the left side of the diagram of
$Q_1$ produces a trivial knot, and the linking number $lk$ of the
corresponding link equals $2+a-b$.
Moreover, the Seifert diagram of $Q_1$ at the bottom of figure~3 is
quasipositive. Here a Seifert diagram is quasipositive if its set of
crossings can be partitioned into single crossings and pairs of crossings, such that the following three conditions are satisfied.
\begin{enumerate}
\item Each single crossing is positive.

\item Each pair of crossings consists of one positive and one negative
crossing joining the same two Seifert circles.

\item A pair of crossings does not separate other pairs of crossings. More
precisely, going from one crossing of a pair to its opposite counterpart along
a Seifert circle, one cannot meet only one crossing of a pair. Such pairs of
crossings are called \emph{conjugating pairs of crossings}.
\end{enumerate}
In \cite{Ba}, we proved that a quasipositive diagram represents a quasipositive
knot. Hence $Q_1$ is a quasipositive knot.

In the second step, we arrange the Vassiliev invariants of order three. Since
'all' the Vassiliev invariants of order less than or equal to two of $Q_1$ and
$K$ coincide (i.e. $v_2(Q_1)=v_2(K)$), we conclude that $Q_1$ and $K$ are
related by a sequence of $C_3$-moves. This is K.~Habiro's result for $n=2$
(see \cite{Ha}). Let $K_1=Q_1$, $K_2$, $\ldots$, $K_l=K$ be a sequence of
knots, such that two succeedding knots are related by a $C_3$-move. Our aim is
to replace this sequence of knots by a sequence of quasipositive knots
$\widetilde K_1=Q_1$, $\widetilde K_2$, $\ldots$, $\widetilde K_l$, such that
$$v_3(\widetilde K_{i+1})-v_3(\widetilde K_i)=v_3(K_{i+1})-v_3(K_i),$$
$1 \leqslant i \leqslant l-1$. By Ohyama and Tsukamoto's result, 
$|v_3(K_2)-v_3(Q_1)|$ depends only
on the permutation $\sigma \in S_{3}$ defined by the cyclic order of the four
strands of the section where the $C_3$-move takes place, as explained above.
From this viewpoint, i.e. if we are only interested in the change of the
Vassiliev invariants of order three, there are only finitely many
combinatorial patterns of $C_3$-moves. A 'standard' pattern of a $C_3$-move
can be applied inside a local box on the right side of the diagram of $Q_1$,
as shown in figure~4. 

\begin{figure}[ht]
\pic{-0pt}{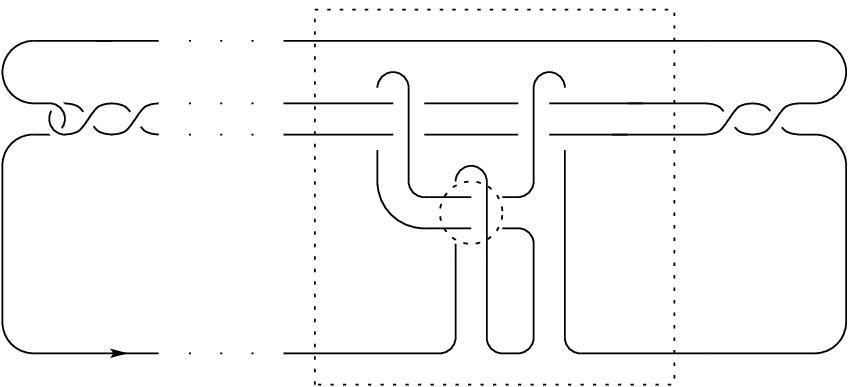} 
\caption{}
\end{figure}

Moreover, we can choose a quasipositive representative for this pattern, i.e.
a representative whose Seifert diagram (inside the box) is quasipositive, see
figure~5. Here we remark that the two segments above and below the
cross-shaped Seifert circle belong to the same Seifert circle since they are
connected outside the local box.

\begin{figure}[ht]
\pic{-0pt}{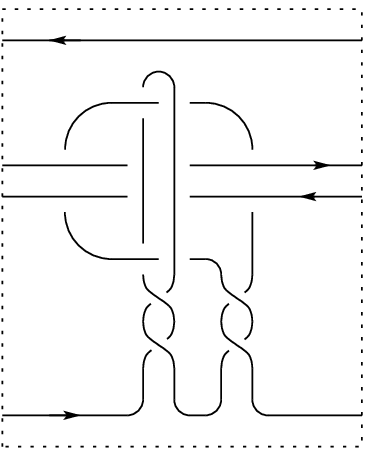} \qquad \qquad \pic{-0pt}{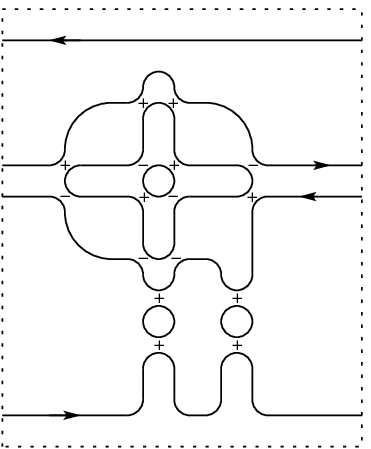}
\caption{}
\end{figure}

However, this standard pattern corresponds to one specific permutation $\sigma
\in S_{3}$. In order to get patterns corresponding to other permutations, we
have to permute the strands inside the local box, as shown by two examples on
the left side of figure~6. We observe that all these patterns have
quasipositive representatives. They are depicted on the right side of
figure~6, together with their Seifert diagrams.

\begin{figure}[ht]
\smallpic{-0pt}{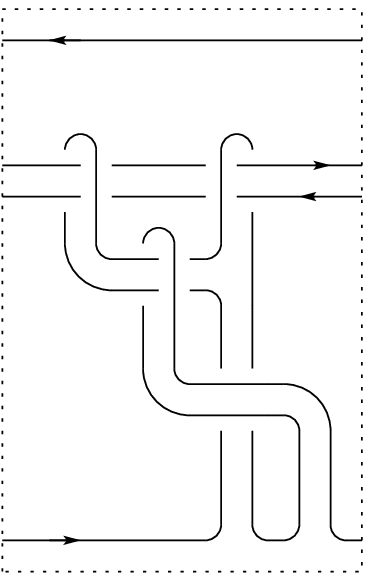} \qquad \smallpic{-0pt}{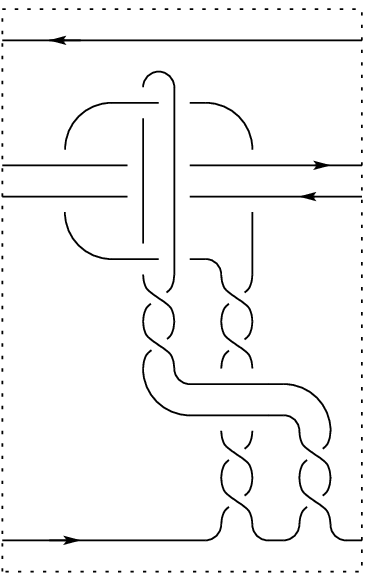} \qquad
\smallpic{-0pt}{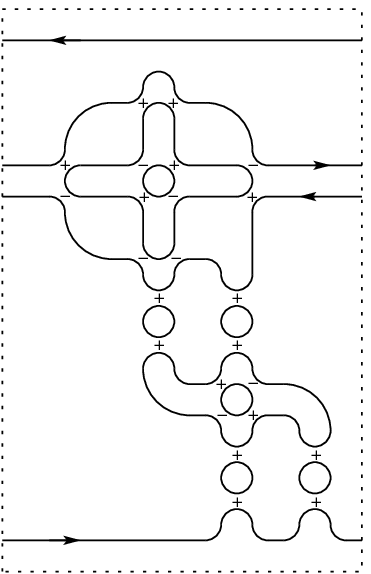}

\bigskip
\bigskip
\smallpic{-0pt}{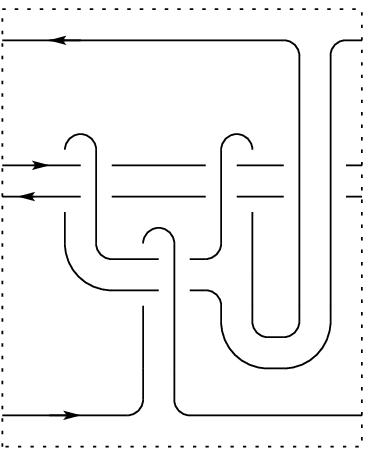} \qquad \smallpic{-0pt}{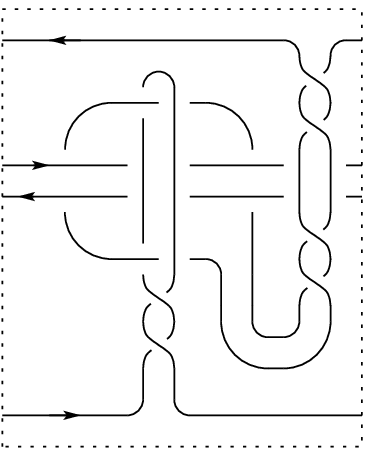} \qquad
\smallpic{-0pt}{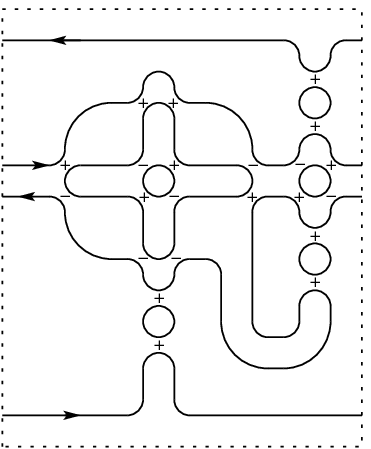}
\caption{}
\end{figure}

Thus we can replace the knot $K_2$ by a quasipositive knot $\widetilde K_2$,
such that
$$v_3(\widetilde K_2)-v_3(Q_1)=\pm(v_3(K_2)-v_3(Q_1)).$$
If the sign of this difference is wrong (i.e. '$-$'), we may arrange it to be
'$+$' by changing four crossings between two strands inside the local box, see
figure~4, where the four crossings are encircled. This modified pattern has
the inverse effect on Vassiliev invariants of order three, as follows from
Ohyama and Tsukamoto's calculation (\cite{OT}, proof of theorem~1.2).

Likewise, we can replace all $C_3$-moves of the sequence $K_1=Q_1$,
$K_2$, $\ldots$, $K_l=K$ by $C_3$-moves that take place in a clearly arranged
box and preserve the quasipositivity of the knot $Q_1$. In this way, we obtain
a sequence of quasipositive knots
$\widetilde K_1=Q_1$, $\widetilde K_2$, $\ldots$, $\widetilde K_l$
and end up with a quasipositive knot
$Q_2:=\widetilde K_l$ whose Vassiliev invariants of order two and three
coincide with those of $K$.

At this point we merely sketch how the process continues: in the $i$-th step,
we arrange the Vassliev invariants of order $i+1$ and define a quasipositive
knot $Q_i$ whose Vassiliev invariants of order less than or equal to $i+1$
coincide with those of $K$. For this purpose, we need only observe that every
combinatorial pattern of a $C_{i+1}$-move has a quasipositive representative:
the heart of such a representative consists of $i^2+i$ conjugating pairs of
arcs. This is illustrated for a $C_4$-move in figure~7.

\begin{figure}[ht]
\pic{-0pt}{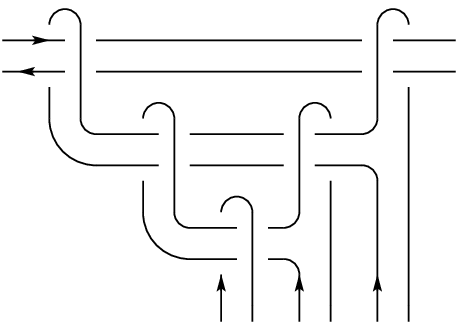} \qquad $\cong$ \qquad \pic{-0pt}{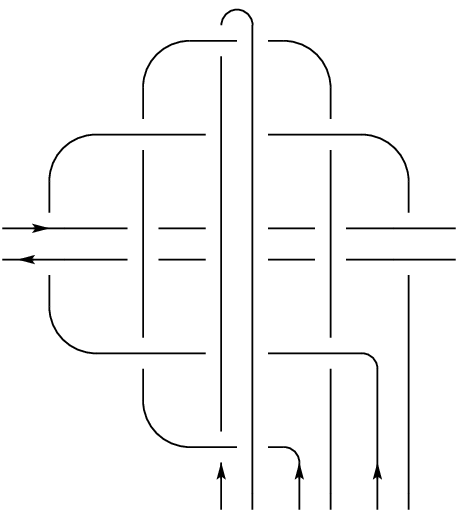}
\caption{}
\end{figure}

At last, the quasipositive knot $Q:=Q_{n-1}$ has the required properties.
\end{proof}

\vbox{%
\begin{remarks}\quad
\begin{enumerate}
\item All the quasipositive knots $Q_i$ can be unknotted by a single crossing
change at the clasp that appears on the left side of their defining diagram
(see e.g. figure~3). In particular, the unknotting number of $Q$ is one, unless $Q$ happens to be the trivial knot.

\item By theorem~1 and Habiro's result, we conclude that every knot can be
transformed into a quasipositive knot by a finite sequence of $C_n$-moves, for
any fixed natural number $n$. It would be interesting to have a direct proof
for this fact, which in turn implies theorem~1. This would possibly simplify
the construction of the desired quasipositive knots.

\item The knot $Q$ might even be chosen to be strongly quasipositive. However,
we do not know how to prove that.
\end{enumerate}
\end{remarks}
}

\textbf{Acknowledgements}. Alexander Stoimenow pointed out an erroneous
statement in the first version of this paper. I would like to thank him for
that and for other useful remarks.

\bigskip
\noindent
Department of Mathematics, University of Basel, Rheinsprung 21, CH-4051 Basel,
Switzerland

\noindent
\emph{E-mail address}: baader@math-lab.unibas.ch


\bigskip
\begin{thebibliography}{99}
   \bibitem{Ba}
     S.~Baader: \emph{Slice and Gordian numbers of track knots}, Osaka J.
     Math. $\mathbf{42}$ (2005).
     
   \smallskip 
   \bibitem{BO}
     M.~Boileau, S.~Orevkov: \emph{Quasi-positivit\'e d'une courbe analytique
     dans une boule pseudo-convexe}, C.~R.~Acad. Sci. Paris S\'er.~I Math
     $\mathbf{332}$ (2001), no.~9, 825-830.

   \smallskip
   \bibitem{Ha}
     K.~Habiro: Master thesis, University of Tokyo (1994).
        
   \smallskip
   \bibitem{O}
     Y.~Ohyama: \emph{Web diagrams and realization of Vassiliev
     invariants by  knots}, J.~Knot Theory Ramifications~$\mathbf{9}$
     (2000), no.~5, 693-701.

   \smallskip
   \bibitem{OTY}
     Y.~Ohyama, K.~Taniyama, S.~Yamada: \emph{Realization of Vassiliev
     invariants by unknotting number one knots}, Tokyo J.~Math~$\mathbf{25}$
     (2002), no.~1, 17-31.
 
   \smallskip
   \bibitem{OT}
     Y.~Ohyama, T.~Tsukamoto: \emph{On Habiro's $C_n$-moves and Vassiliev
     invariants of order $n$}, J.~Knot Theory Ramifications~$\mathbf{8}$
     (1999), no.~2, 15-26.
     
   \smallskip
   \bibitem{Ru1}
     L.~Rudolph: \emph{Constructions of quasipositive knots and links. I}.
     Knots, braids and singularities (Plans-sur-Bex, 1982), 233-245, Monogr.
     Enseign. Math. $\mathbf{31}$, Enseignement Math. (Geneva, 1983).

   \smallskip
   \bibitem{Ru2}
     L.~Rudolph: \emph{Algebraic functions and closed braids}, Topology
     $\mathbf{22}$ (1983), no.~2, 191-202.

   \smallskip
   \bibitem{St}
     A.~Stoimenow: \emph{Vassiliev invariants and rational knots of unknotting
     number one}, Topology $\mathbf{42}$ (2003), no.~1, 227-241.
\end{thebibliography}
\end{document}